\documentclass[12pt]{article}
\usepackage{amsmath} 
\usepackage{epsfig}
\usepackage{amssymb} 
\newtheorem{exercise}{Exercise}
\newtheorem{solution}{Exercise}
\newtheorem{lm}{{\bf Lemma}}[section]
\newtheorem{theor}[lm]
{{\bf Theorem}}
\newtheorem{deff}[lm]{{\bf Definition}}
\newtheorem{cor}[lm]
{{\bf Corollary}}
\newtheorem{prop}[lm]{{\bf Proposition}}
\newtheorem{remark}[lm]{Remark}

\newcommand{\lem}[1]{\begin{lm}$\!\!$\sep{\sl #1}\end{lm}}
\newcommand{\theo}[1]{\begin{theor}$\!\!$\sep{\sl #1}\end{theor}}
\newcommand{\coro}[1]{\begin{cor}$\!\!$\sep{\sl #1}\end{cor}}

\newcommand{\propo}[1]{\begin{prop}$\!\!$\sep{\sl #1}\end{prop}}
\newcommand{\rem}[1]{{\begin{remark}$\!\!$\sep{\rm #1}\end{remark}}}

\newcommand{\rf}[1]{{\rm(\ref{#1})}}

\newcommand{\pr}{\noindent {\sc Proof}. --- \ }

\newcommand{\ep}{\hfill \framebox[2mm]{\ } \medskip}
\newcommand{\sep}{{\rm .} --- \ }
 
\newcommand{\be}{\begin{enumerate}}
\newcommand{\ee}{\end{enumerate}}

\newcommand{\ds}{\displaystyle}

\newcommand{\N}{\mathbb{N}}

\newcommand{\R}{\mathbb{R}}

\newcommand{\g}{{\bf g}}

\begin{document} 

\title{\bf On the solutions of a boundary value problem arising in free convection with prescribed heat flux}
\author{\sc Mohamed A\"IBOUDI and Bernard BRIGHI\footnote{~Corresponding author.}}
\date{}
\maketitle
\bigskip

\begin{abstract}
\noindent For given $a\in\R$, $c<0$, we are concerned with the solution $f^{}_b$ of the differential equation $f^{\prime\prime\prime}+ff^{\prime\prime}+\g(f^{\prime})=0$, satisfying the initial conditions $f(0)=a$, $f'(0)=b$, $f''(0)=c< 0$,  where $\g$ is some nonnegative subquadratic locally Lipschitz function.
It is proven that there exists $b_*>0$ such that $f^{}_b$ exists on $[0,+\infty)$ and is such that $f'_b(t)\to 0$ as $t\to+\infty$, if and only if $b\geq b_*$. This allows to answer questions about existence, uniqueness and boundedness of solutions to a boundary value problem arising in fluid mechanics, and especially in boundary layer theory.


\end{abstract}

\section{Introduction.}
\noindent
We consider the similarity third order differential equation 
\begin{equation}
f'''+ff''+\g(f')=0\label{eqg}
\end{equation}
on $[0,+\infty)$, with the boundary conditions
\begin{align}
&f(0)=a, \label{c1} \\
&f''(0)=c<0, \label{c2} \\
&f'(+\infty):=\ds\lim_{t\to+\infty}f'(t)=0, \label{c3}
\end{align}
where the function $\g:\R\to\R$ is assumed to be  locally Lipschitz.

This boundary value problem with $\g(x)=\beta x^2$ arises in fluid mechanics, when looking for similarity solutions in free convection boundary-layer flows adjacent to permeable surfaces in porous media.
The initial condition \rf{c2} means that heat flux is prescribed on the surface. 
In other situations, the surface temperature is prescribed, and in this case, condition \rf{c2} has to be replaced by 
$f'(0)=b>0$.
See for example \cite{pop1} and \cite{pop} for details on the derivation of these problems, in the context of the boundary layer theory.

For $\g(x)=\beta x^2$, mathematical analysis of the problem with prescribed surface temperature is done in \cite{brighi02}, \cite{brighi01}, \cite{brighisari}, \cite{guedda1}, \cite{guotsai} and \cite{tsai}. See also \cite{amuc} and \cite{bb} for general function $\g$. With prescribed surface heat flux, see \cite{heat_flux} and \cite{tsaiwang}. 

In this paper, we are interested in the boundary value problem \rf{eqg}-\rf{c3}, with $0\leq\g(x)\leq x^2$. The particular  case where $\g(x)=\beta x^2$ with $0<\beta<1$ corresponds to a question, which was not solved in \cite{heat_flux}, and which has obtained an answer in \cite{tsaiwang}. The method used by J.-C. Tsai and C.-A. Wang is based on the fact that $\g$ is homogeneous of degree 2, and on the study of a plane vector field associated to the differential equation \rf{eqg}. Here, we propose to revisit this question in a direct way, and, as far as possible, to prove results with $\g$ such that $0\leq\g(x)\leq x^2$. We will see that, under this hypothesis, we are able to get existence of solutions, but that we will need to assume that $\g(x)=\beta x^2$ with $0<\beta<1$ to get more precise results (as uniqueness of the bounded solution).
However, we think that this latter assumption is not necessary. In fact, for the boundary value problem involving \rf{eqg} and the boundary conditions corresponding to prescribed surface temperature, it is possible to prove that the bounded solution is unique, only by assuming that $0\leq\g(x)\leq x^2$, see \cite{bb}.

\section{Preliminary remarks.}

\noindent 
The method to solve the boundary value problem \rf{eqg}-\rf{c3} is shooting. For that,
let  $f^{}_b$ denote the solution of the initial value problem :\begin{equation*}
\left \{ \begin{array}
[c]{lll}
f'''+ff''+\g(f')=0,\\ \noalign{\vskip1mm}
f(0) = a,\\ \noalign{\vskip1mm}
f'(0) =b,\\ \noalign{\vskip1mm}
f''(0) = c,
\end{array} \right.\leqno ({\cal P}_{\g;a,b,c})
\end{equation*}
and let $[0,T_b)$ be the right maximal interval of existence of $f^{}_b$. To obtain a solution of the boundary value problem \rf{eqg}-\rf{c3} amounts to find a value of $b$ such that $T_b=+\infty$ and $f'_b(t)\to0$ as $t\to+\infty$.

The following useful identities are obtained, by multiplying equation \rf{eqg} by $1$, $f^{}_b$ and $t$ respectively, and integrating by parts. For all $t\in[0,T_b)$, we have :
\begin{equation}
f''_b(t)-c+f^{}_b(t)f'_b(t)-ab=\int_0^t \left(f'_b(s)^2-\g(f'_b(s))\right)ds \label{i1}
\end{equation}
\begin{equation}
f^{}_b(t)f''_b(t)-ac-\frac12f'_b(t)^2+\frac12 b^2+f^{}_b(t)^2f'_b(t)-a^2b=\int_0^t f^{}_b(s)\left(2f'_b(s)^2-\g(f'_b(s))\right)ds \label{i2}
\end{equation}
and
\begin{equation}
tf''_b(t)-f'_b(t)+b+tf^{}_b(t)f'_b(t)-\frac12f^{}_b(t)^2+\frac12 a^2=\int_0^t s\left (f'_b(s)^2-\g(f'_b(s))\right)ds. \label{i3}
\end{equation}

\medskip
We now give some lemmas, that we will use in the next sections. 
\bigskip
\lem{\label{Tfini} If $T_b<+\infty$, then $f''_b(t)$ and $f'_b(t)$ are unbounded as $t\to T_b$.
}
\pr
First, if $T_b$ is finite, then $\vert f^{}_b(t)\vert+\vert f'_b(t)\vert+\vert f''_b(t)\vert$ is unbounded as $t\to T_b$. Then, necessarily $\vert f''_b(t)\vert$ is unbounded as $t\to T_b$, and using \rf{i1} we see that $\vert f'_b(t)\vert$ is also unbounded as $t\to T_b$.
\ep
\rem{In general, $f^{}_b$ has no reason to be unbounded ; for example
$f:t\mapsto\sqrt{1-t}$ is the solution of $({\cal P}_{\g;1,-1/2,-1/4})$, with $\g(x)=x^2\left(1-12x^3\right)$, on the maximal interval $[0,1)$ and $f''(t)\to-\infty$, $f'(t)\to-\infty$ and 
$f(t)\to 0$ as $t\to 1$.
}

\medskip
\lem{\label{concave} If $\g(x)>0$ for $x\not=0$ and if $c<0$, then for any $b\in\R$ we have $f''_b<0$ on $[0,T_b)$, {\sl i.e.} $f^{}_b$ is concave.
}
\pr This follows from the relation 
$
\left(f''e^F\right)'=-\g(f')e^F\label{exp}
$, where $F$ is any primitive function of $f$ on $[0,T_b)$.
\ep
%

%
%

%
\rem{\label{blaineq} It is possible to show, under the assumptions of the previous lemma, that if $T_b=+\infty$ then $f'_b>0$. Indeed, on the contrary, there would exists $t_0>0$ such that $f^{}_b$ and $f'_b$ are negative on $(t_0,+\infty)$ and therefore $f^{}_b$ would be a negative concave subsolution of the Blasius equation ({\sl i.e.} satisfying $f'''+ff''\leq 0$) on $(t_0,+\infty)$, and using similar arguments as the ones in the proofs of Proposition 2.1 and 2.2 of \cite{brighi03}, we would obtain a contradiction. See also \cite{bb}.
}


\medskip

\lem{\label{a-} Let us assume that $0<\g(x)\leq 2x^2$ for $x\not=0$, $a<0$ and $c<0$. If  $b>0$ is large enough, then there exists $s_b\in(0,T_b)$ such that $f^{}_b(s_b)=0$ and $f'_b(s_b)>\frac{3b}{4}$. 
}
\pr Since $f''_b(0)=c<0$, we deduce from Lemma \ref{concave} that $f''_b<0$ on $[0,T_b)$.
Let us assume that there exists $t\in(0,T_b)$ such that $f^{}_b(t)\leq0$ and $f'_b(t)=\frac{3b}{4}$. Then, using \rf{i2}, we have 
$$-ac+\frac{7b^2}{32}-a^2b=-f^{}_b(t)f''_b(t)-f^{}_b(t)^2f'_b(t)+\int_0^t f^{}_b(s)\left(2f'_b(s)^2-\g(f'_b(s))\right)ds\leq 0$$ and then $b$ is smaller than the positive root of the polynomial $
7X^2-32 a^2X-32ac$.
This completes the proof.
\ep
\section{The solutions of $({\cal P}_{\g;a,b,c})$ when $\g$ is nonnegative and subquadratic.}

\noindent In this section, we will assume that $c<0$ and that $\g:\R\to\R$ is locally Lipschitz and such that $0<\g(x)\leq x^2$ for all $x\not=0$. Let us notice that, by continuity, we have $\g(0)=0$.

By Lemma \ref{concave},  the function $f^{}_b$ is concave on $[0,T_b)$, for all $b\in\R$. We will distinguish the following two types of behavior.
\begin{itemize}
\item {\sl Type} (I) : $f'_b\geq 0$ on $[0,T_b)$.

\item {\sl Type} (II) : there exists $t_0\in[0,T_b)$ such that $f'_b<0$ on $(t_0,T_b)$.
\end{itemize}
We then define the sets 
$${\cal B}_1=\{b\in\R~;~f^{}_b\mbox{ is of type (I)}\}~~~\mbox{ and }~~~{\cal B}_2=\{b\in\R~;~f^{}_b\mbox{ is of type (II)}\}.$$
Clearly, we have ${\cal B}_1\cap{\cal B}_2=\emptyset$ and ${\cal B}_1\cup{\cal B}_2=\R$. Moreover, ${\cal B}_1$ is a closed set (and hence ${\cal B}_2$ is an open set). In fact, if $b_n\in{\cal B}_1$ is a sequence converging to some $b_*\in\R$, and if $t\in[0,T_{b_*})$, then, from the lower semicontinuity of the mapping $b\mapsto T_b$, there exists $n_0\in\N$ such that for $n\geq n_0$ we have $T_{b_n}>t$. Now, 
the continuity of  $(b,t)\mapsto f'_b(t)$ (defined for $b\in\R$ and $t\in[0,T_b)$), allows to write 
$$f_{b_*}(t)=\lim_{n\to+\infty}f_{b_n}(t)\geq 0,$$  
and thus $b_*\in{\cal B}_1$.

On the other hand, it is clear that $(-\infty, 0]\subset{\cal B}_2$, or equivalently that ${\cal B}_1\subset(0,+\infty)$. {A priori}, nothing indicates that ${\cal B}_1\not=\emptyset$.

\rem{It follows from Remark \ref{blaineq}, that, if $b\in{\cal B}_2$, then $T_b<+\infty$.
}

\medskip
The following result gives informations about $f^{}_b$ for $b\in{\cal B}_1$.

\propo{\label{mu+} 
If $b\in{\cal B}_1$, then $T_b=+\infty$, $f'_b>0$, $f'_b(t)\to 0$ as $t\to+\infty$, and there exists $t_0\geq 0$ such that $f^{}_b(t)>0$ for all $t>t_0$. If, in addition, $f^{}_b$ is bounded, then there exists a positive constant $A_b$ such that : 
\begin{align}
\label{simf}&f^{}_b(t)=\mu_b-A_be^{-\mu_b t(1+{\rm o}(1))}\\ \noalign{\vskip2mm}
\label{simff}&f'_b(t)\sim \mu_b(\mu_b-f^{}_b(t))~~~\mbox{ and }~~~f''_b(t)\sim-\mu_b f'_b(t)
\end{align}
as $t\to+\infty$, where $\mu_b>0$ is the limit of $f^{}_b(t)$ as $t\to+\infty$. 
}

\pr 
Let $b\in{\cal B}_1$. Since $f''_b<0$ on $[0,T_b)$, then $f'_b(t)$ has a nonnegative limit $\ell$ as $t\to T_b$ and thanks to Lemma \ref{Tfini}, it follows that $T_b=+\infty$. 
Now, we claim that $\ell=0$\footnote{~See Lemma 3 of \cite{amuc} for a general proof of the fact that $\g(\ell)=0$.}. In fact, if $\ell>0$, then we have $f'_b(t)\sim\ell$ and $f^{}_b(t) \sim \ell t$ as $t \to +\infty$. Using \rf{i1}, since $\g(l)>0$, we get
\begin{align*}
f''_b(t)-c-ab&=-f^{}_b(t)f'_b(t)+\int_{0}^t f'_b(s)^2ds-\int_{0}^t\g(f'_b(s))ds \cr
&=-\ell^2t(1+{\rm o}(1))+\ell^2t(1+{\rm o}(1))-\g(\ell)t(1+{\rm o}(1)) \cr
&=-\g(\ell)t+{\rm o}(t)
\end{align*}
which implies that $f''_b(t)\to-\infty$ as $t\to +\infty$. This is a contradiction with the fact that $f'_b(t)$ has a finite limit as $t\to +\infty$. Therefore, $\ell=0$, and the positivity of $f'_b$ then follows from the fact that $f''_b<0$.

Now, if $a\geq 0$ then $f^{}_b$ is positive, and if $a<0$ and $f^{}_b(t)<0$ for all $t\geq 0$ then we have
$$\forall t\geq0,~~~f'''_b=-f^{}_bf''_b-\g(f'_b)<0$$
whence we deduce that $f'_b$ is concave. This contradicts the fact that $f'_b$ is positive and tends to $0$ at infinity. Therefore, there exists $t_0\geq 0$ such that $f^{}_b(t)>0$ for all $t>t_0$.

\smallskip
Finally, if $f^{}_b$ is bounded, then we can integrate \rf{eqg} between $t$ and $+\infty$ ; this gives :
$$f''_b(t)+f^{}_b(t)f'_b(t)=\int_t^{+\infty}\left(\g(f'_b(s))-f'_b(s)^2\right)ds$$ 
Together with the assumptions about $\g$ and the fact that $f'_b$ is positive and decreasing, we get :
\begin{align*}
0\geq\frac{f''_b(t)}{f'_b(t)}+f^{}_b(t)&\geq-\frac{1}{f'_b(t)}\int_t^{+\infty}f'_b(s)^2ds
\geq-\int_t^{+\infty}f'_b(s)ds=f^{}_b(t)-\mu_b. 
\end{align*}
We immediatly deduce the second relation of \rf{simff} ; the first one  follows from the L'H\^opital's rule, and \rf{simf} by an integration.
\ep

\medskip
\propo{\label{b1nonvide} 
Let $\delta\in(0,1)$ and let us assume that $0<\g(x)\leq (1-\delta)x^2$ for all $x\not=0$
. There exists $b_0>0$ such that $[b_0,+\infty)\subset{\cal B}_1$.
}
\pr As previously, let $f^{}_b$ denote the solution of the initial value problem $({\cal P}_{\g;a,b,c})$. 
Taking into account Lemma \ref{a-} and Proposition \ref{mu+}, we deduce that, for $b$ large enough, there exists $t_b\in(0,T_b)$ such that $f'(t_b)=\frac b2$ and $f^{}_b(t_b)\geq 0$. Using \rf{i1} and \rf{i3}, we can write :
\begin{align*}
t_bf''_b(t_b)+t_bf^{}_b(t_b)f'_b(t_b)-(c+ab)t_b&=\int_0^{t_b} t_b\left(f'_b(s)^2-\g(f'_b(s))\right)ds, \\ \noalign{\vskip2mm}
t_bf''_b(t_b)+\frac b2+t_bf^{}_b(t_b)f'_b(t_b)-\frac12f^{}_b(t_b)^2+\frac12 a^2&=\int_0^{t_b} s\left (f'_b(s)^2-\g(f'_b(s))\right)ds.
\end{align*}
By substraction, we obtain
$$-(c+ab)t_b-\frac b2+\frac12f^{}_b(t_b)^2-\frac12 a^2=\int_0^{t_b} (t_b-s)\left (f'_b(s)^2-\g(f'_b(s))\right)ds\geq 0.$$
Whence, together with the fact that $0\leq f^{}_b(t_b)\leq bt_b+a$, we get
$$a^2+b+2(c+ab)t_b\leq f^{}_b(t_b)^2\leq (bt_b+a)^2$$
and hence $b^2t_b^2-2ct_b-b\geq 0$. It follows that 
$$t_b\geq \frac{c+\sqrt{c^2+b^3}}{b^2}
=\frac{b}{\sqrt{c^2+b^3}-c}.$$
Then, if $f'_b(t)=0$, we have $t>t_b$ and, using \rf{i1} and the assumptions about $\g$, we have
\begin{align*}
0>f''_b(t)&=c+ab+\int_0^t \left(f'_b(s)^2-\g(f'_b(s))\right)ds \\ \noalign{\vskip2mm}
&\geq c+ab+\delta\int_0^t f'_b(s)^2ds\geq c+ab+\delta\int_0^{t_b} f'_b(s)^2ds\\ \noalign{\vskip2mm}
&\geq c+ab+\frac \delta4 t_b b^2\geq c+ab+\frac \delta4 \frac{b^3}{\sqrt{c^2+b^3}-c}.
\end{align*}
It follows that, if $b$ is large enough, then $f'_b$ does not vanish, and hence $b\in{\cal B}_1$.
\ep

\rem{If we only suppose that $\g(x)\leq x^2$, then the previous result may not hold ; for example, if $\g(x)=x^2$ and $a\leq 0$, $c<0$, then ${\cal B}_1=\emptyset$. Indeed, if there exists $b\in{\cal B}_1$, then $b>0$, and identity \rf{i1} gives $f''_b(t)+f^{}_b(t)f'_b(t)=c+ab$ for all $t\in\R$.
Integrating, we obtain :
$$\forall t\geq 0,~~~f'_b(t)+\frac12f^{}_b(t)^2=b+\frac12 a^2+(c+ab)t$$
and we see that the right hand side tends to $-\infty$ as $t\to+\infty$, whereas the left one is positive. We have a contradiction.}

\medskip

\propo{\label{borne} If $b$ is a point of the boundary of ${\cal B}_1$, then $f^{}_b$ is bounded. 
}
\pr
\noindent Since $b$ is on the boundary of ${\cal B}_1=\R\setminus{\cal B}_2$, there exists a sequence of positive real numbers $b_n\in{\cal B}_2$ converging to $b$. Let us set $f_n=f_{b_n}$ and $T_n=T_{b_n}$. Since $b_n\in{\cal B}_2$, there exists $t_n\in(0,T_n)$ such that $f'_n(t_n)=0$.

First, we remark that $t_n\to+\infty$ as $n\to+\infty$. On the contrary, there would exist an increasing subsequence $t_{n_k}$ converging to some $t<+\infty$. 
Because of the lower semicontinuity of the mapping $b\mapsto T_b$,
for $k$ large enough, we should have $T_{n_k}>t$ and we could write 
$$0=\lim_{k\to+\infty}f'_{n_k}(t_{n_k})\geq \lim_{k\to+\infty}f'_{n_k}(t)=f'_b(t)$$
and get a contradiction.

Next, using \rf{i3} for $t=t_n$ yields
$$b-\frac12\left(f^{}_n(t_n)^2-a^2\right)=-t_nf''_n(t)+\int_0^{t_n} s\left(f'_n(s)^2-\g(f'_n(s)\right)ds\geq0$$
and hence 
$f^{}_n(t_n)^2\leq 2b+a^2.$
Since $f'_n$ is positive on $[0,t_n)$, we obtain
$$\forall s\in[0,t_n],~~~f^{}_n(s)\leq \sqrt{2b+a^2}.$$
To conclude, we fix $t\in[0,+\infty)$. For $n$ large enough, we have $t_n>t$ and 
$$f^{}_b(t)=\lim_{n\to+\infty}f_{n}(t)\leq\sqrt{2b+a^2}.$$
This completes the proof.
\ep

\section{The boundary value problem \rf{eqg}-\rf{c3} when $\g(x)=\beta x^2$ with $0<\beta<1$.}

\noindent 
Here, we consider the case where  $\g(x)=\beta x^2$ with $0<\beta<1$. Our main result is the following.

\theo{\label{beta1} Let $a\in\R$ and $c<0$. If $\g(x)=\beta x^2$ with $0<\beta<1$, then there exists $b_*>0$ such that ${\cal B}_1=[b_*,+\infty)$. Moreover, if $b>b_*$, then $f^{}_b$ is unbounded.
}
\pr Taking into account Propositions \ref{b1nonvide} and \ref{borne}, it is sufficient to prove that there is at most one $b>0$ such that $f^{}_b$ is bounded.

First, let us assume that for some $b>0$, the function $f=f^{}_b$ is bounded. Let $\mu>0$ be the limit of $f(t)$ as $t\to+\infty$. Since $f$ is concave and increasing, then we can define a function $v:(0,1]\rightarrow\R$, such that
$$\forall t\geq 0,~~~v\left(\frac{1}{b^2}f'(t)^2\right)=\frac1{\sqrt b}f(t).$$
By setting $y=\ds\frac{1}{b^2}f'(t)^2$, we get 
\begin{equation}
f(t)=\sqrt{b}\,v(y),~~~f'(t)=b\sqrt{y},~~~f''(t)=\frac{b^{3/2}}{2v'(y)}~~\text{ and }~~~f'''(t)=-\frac{b^2v''(y)\sqrt{y}}{2v'(y)^3}.\label{y1}
\end{equation}
Then, using \rf{eqg} we obtain
\begin{equation}
\forall y\in(0,1],~~~v''(y)=\frac{v(y)v'(y)^2}{\sqrt{y}}+2\beta\sqrt{y}\,v'(y)^3. \label{eq-v1}
\end{equation}
Moreover $v$ is decreasing on $(0,1]$, and we have
\begin{equation}
v(1)=\frac{a}{\sqrt b},~~~v'(1)=\frac{b^{3/2}}{2c}~~\mbox{ and }~~\lim_{y\to 0}v(y)=\frac{\mu}{\sqrt b}.\label{init1}
\end{equation}
In addition, using \rf{simff}, it holds
\begin{equation}
v'(y)\sim-\frac{\sqrt b}{2\mu\sqrt y}~~\mbox{ as }~~y\to 0.\label{init2}
\end{equation}

Now, let us assume that there are $b_1>b_2$ such that the functions $f_1=f^{}_{b_1}$ and $f_2=f^{}_{b_2}$ are bounded. For $i=1,2$, let $\mu_i$ be the limit of $f_i(t)$ as $t\to+\infty$. Let $v_i:(0,1]\to\R$ be the corresponding solutions of  \rf{eq-v1}. 

If $w=v_1-v_2$, then $w$ is defined on $(0,1]$ and we have
\begin{equation}
w(1)=\frac{a}{\sqrt b_1}-\frac{a}{\sqrt b_2},~~~w'(1)=\frac1{2c}\left(b_1^{3/2}-b_2^{3/2}\right)<0~~\mbox{ and }~~\lim_{y\to 0}w(y)=\frac{\mu_1}{\sqrt b_1}-\frac{\mu_2}{\sqrt b_2}.\label{b2}
\end{equation}
Moreover, $w$ cannot have neither positive maximum, or negative minimum in $(0,1)$. Indeed, if $x\in(0,1)$ is such that $w(x)>0$, $w'(x)=0$, then using \rf{eq-v1} we have
$$w''(x)=\frac{v'_1(x)^2}{\sqrt{x}}w(x)>0.$$ 
The same arguments show that $w$ has no negative minimum in $(0,1)$.

We now distinguish between the cases $a\leq 0$ and $a>0$. 
\begin{itemize}
\item If $a\leq 0$, then $w(1)\geq 0$. Since $w$ cannot have a positive maximum in $(0,1)$ and $w'(1)<0$, it follows that $w(0)>0$. Thus,  $\frac{\mu_1}{\sqrt b_1}>\frac{\mu_2}{\sqrt b_2}$ and hence, using \rf{init2}, we get 
\begin{equation*}
w'(y)\sim-\frac1{2\sqrt y}\left(\frac{\sqrt b_1}{\mu_1}-\frac{\sqrt b_2}{\mu_2}\right)~~\mbox{ as }~~y\to 0.\label{ww}
\end{equation*}
Therefore, $w'(y)\to+\infty$ as $y\to 0$ and this gives a contradiction since, on the contrary,  $w$ should have a positive maximum in $(0,1)$.

\item If $a>0$, then $w(1)<0$. Using the fact that $w$ cannot have neither positive maximum, or negative minimum in $(0,1)$ and the same arguments as previously, we obtain that necessarily we have $w'\leq 0$ and $w(1)<w(0)\leq 0$.

Now, for $i=1,2$, let us set $V_i=\frac1{v'_i}$ and  $W=V_1-V_2$.  First, thanks to \rf{init2}, we have $W(y)\to 0$ as $y\to 0$. Next, using \rf{eq-v1}, we obtain
\begin{align*}
\forall y\in (0,1],~~~W'(y)&=-\frac{v''_1(y)}{v'_1(y)^2}+\frac{v''_2(y)}{v'_2(y)^2}=-\frac{w(y)}{\sqrt{y}}-2\beta \sqrt{y} \, w'(y)\\ \noalign{\vskip2mm}
&\leq-\frac{w(y)}{\sqrt{y}}-2w'(y)\sqrt{y}=-2\bigl(\sqrt{y}w(y)\bigr)'.
\end{align*}
Integrating between $0$  and $1$, we  get
$W(1)\leq -2w(1)$. Thus,
$$2c\left(\frac{1}{b_1^{3/2}}-\frac{1}{b_2^{3/2}}\right)\leq -2a\left(\frac{1}{\sqrt{b_1}}-\frac{1}{\sqrt{b_2}}\right).$$
Hence
$$c\left(\frac{1}{b_1}+\frac{1}{\sqrt{b_1b_2}}+\frac{1}{b_2}\right)\geq -a$$
and
$$\frac{c}{b_1}+a\geq-c\left(\frac{1}{\sqrt{b_1b_2}}+\frac{1}{b_2}\right)>0.$$
This is contradiction, since thanks to \rf{i1} written for $f_1$ and $t\to+\infty$ we must have
$c+ab_1<0$. 
\end{itemize}
The proof is complete.\ep

\rem{If $b>b_*$, then there exists a positive constant $A_b$ such that 
$$f_b(t)\sim A_b t^{\frac1{1+\beta}}~~\mbox{ as }~t\to+\infty.$$
See \cite{equiv} and \cite{guedda1}.}

\medskip

\coro{\label{beta} Let $a\in\R$ and $c<0$. If $\g(x)=\beta x^2$ with $0<\beta<1$, then the boundary value problem \rf{eqg}-\rf{c3} has exactly one bounded solution, and infinitely many unbounded solutions.
}
\pr This follows immediatly from the previous theorem.\ep

\medskip
\coro{\label{mm} Let $a\in\R$ and $m\in\left(-1,-\frac12\right)$. The boundary value problem 
\begin{equation}
\left \{ \begin{array}
[c]{lll}
f'''+(m+2)ff''-(2m+1)f'^2=0~\mbox{ on }~[0,+\infty) \\ \noalign{\vskip2mm}
f(0) = a,~~~f''(0) = -1,~~~f'(t)\to 0~\mbox{ as }~t\to+\infty,
\end{array} \right.\label{mma}
\end{equation}
has exactly one bounded solution, and infinitely many unbounded solutions.
}
\pr By setting 
$$\frac1{\sqrt{m+2}}\,\tilde{f}(t)=f \biggl(\frac{t}{\sqrt{m+2}}\biggr)$$
we see that $f$ is a solution of the boundary value problem \rf{mma} if and only if $\tilde{f}$ is a solution of the boundary value problem \rf{eqg}-\rf{c3} with $\g(x)=\beta x^2$ and $\beta=-\frac{2m+1}{m+2}$. The proof then follows from Corollary \ref{beta}.\ep

\section{Conclusion.}

\noindent The result of Theorem \ref{beta1} is obtained by using strongly the fact that $\g$ is homogeneous of degree 2, contrary to the results of sections 2 and 3, where only the subquadratic nature of $\g$ is used. 
Nevertheless, we conjecture that the uniqueness of the bounded solution should hold under the hypothesis $0\leq \g(x)\leq x^2$. In fact, the boundary value problem
\begin{equation*}
\left \{ \begin{array}
[c]{lll}
f'''+ff''+\g(f')=0~\mbox{ on }~[0,+\infty)\\ \noalign{\vskip2mm}
f(0) = a,~~~f'(0) = b>0,~~~f'(t)\to 0~\mbox{ as }~t\to+\infty,
\end{array} \right.
\end{equation*}
has at most one bounded concave solution, if $0\leq \g(x)\leq x^2$, see \cite{bb}, and this result is one of the reasons for which we hope that this conjecture holds.

On the other hand, let us notice that we recover the results of J.-C. Tsai and C.-A. Wang \cite{tsaiwang}, in a totally different way, and perhaps more directly.

\bigskip

\noindent Mohamed A\"iboudi

\noindent Universit\'e d'Oran (Es-Senia)

\noindent D\'epartement de Math\'ematiques

\noindent BP 1524

\noindent El Menouar

\noindent 31000 Oran

\noindent Algeria

\noindent e-mail : {\tt m.aiboudi@yahoo.fr}

\bigskip

\noindent Bernard Brighi

\noindent Universit\'e de HauteAlsace 
 
\noindent Laboratoire de Math\'ematiques, Informatique et Applications
 
\noindent 4 rue des fr\`eres Lumi\`ere 

\noindent F-68093 Mulhouse cedex

\noindent France

\noindent e-mail : {\tt bernard.brighi@uha.fr}


\begin{thebibliography}{99}                                           




\bibitem{brighi03} Z. Belhachmi, B. Brighi and K. Taous, On the concave solutions of the Blasius equation, {\it Acta Math. Univ. Comenian.} {\bf 69} (2000), no. 2, 199-214.


\bibitem {brighi02} Z. Belhachmi, B. Brighi and K. Taous, On a family of differential equations for boundary layer approximations in porous media, {\it European J. Appl. Math.} {\bf 12} (2001), no. 4,
513-528.



\bibitem {brighi01} B. Brighi, On a similarity boundary layer equation, {\it Z. Anal. Anwendungen} {\bf 21} (2002), no. 4, 931-948.

\bibitem {bb} B. Brighi,  On the differential equation $f'''+ff''+\g(f')=0$ and the associated boundary value problems, preprint.

\bibitem {heat_flux} B. Brighi and J.-D. Hoernel, On similarity solutions for boundary layer flows with prescribed heat flux, {\it Math. Methods Appl. Sci.} {\bf 28} (2005), no. 4, 479-503.


\bibitem{equiv} B. Brighi and J.-D. Hoernel, Asymptotic behavior of the unbounded solutions of some boundary layer equation, {\it Arch. Math. (Basel)} {\bf 85} (2005), no. 2, 161-166.



\bibitem {amuc} B. Brighi and J.-D. Hoernel, On a general similarity boundary layer equation, preprint (http://arxiv.org/abs/math/0601385).

\bibitem {brighisari} B. Brighi and T. Sari, Blowing-up coordinates for asimilarity boundary layer equation,  {\it Discrete Contin. Dyn. Syst. (Serie A)} {\bf 12} (2005), no. 5, 929-948.

\bibitem {pop1} M. A. Chaudhary, J.H. Merkin and I. Pop, Similarity solutions in
free convection boundary-layer flows adjacent to vertical permeable surfaces
in porous media: I prescribed surface temperature, {\it European J. Mech. B Fluids} {\bf 14} (1995), no. 2, 217-237.

\bibitem {pop}M. A. Chaudhary, J.H. Merkin and I. Pop, Similarity solutions in
free convection boundary-layer flows adjacent to vertical permeable surfaces
in porous media: II prescribed surface heat flux, {\it Heat and Mass Transfer} {\bf 30}
 (1995), 341-347.








\bibitem{guedda1} M. Guedda, Similarity solutions of differential equations for boundary layer approximations in porous media, {\it J. Appl. Math. Phys.  (ZAMP)} {\bf 56} (2005), 749-762.




\bibitem{guotsai} J.-S. Guo and J.-C. Tsai, The structure of the solutions for a third order differential equation in boundary layer theory,  {\it Japan J. Indust. Appl. Math.} {\bf 22} (2005), no. 3, 311-351.


\bibitem {hart}P. Hartmann, {\it Ordinary Differential Equations}. Wiley, New-York (1964).





















\bibitem{tsai} J.-C. Tsai, Similarity solutions for boundary layer flows with prescribed surface temperature, to appear in {\it Appl. Math. Lett.}, available online 20 March 2007.

\bibitem{tsaiwang} J.-C. Tsai and C.-A. Wang, A note on similarity solutions for boundary layer flows with prescribed heat flux, to appear in {\it Math. Methods Appl. Sci.}, available online 15 March 2007.









\end{thebibliography}
\end{document}